\documentclass[times, 10pt]{article}

\newtheorem{Theorem}{Theorem}

\newtheorem {Lemma} [Theorem]    {Lemma}
\newtheorem {Corollary}  [Theorem]    {Corollary}
\newtheorem {Proposition}[Theorem]    {Proposition}
\newtheorem{Conjecture}[Theorem]{Conjecture}
\newcommand{\proof}{\noindent {\bf Proof.} }

\newcommand{\s}{\sigma}
\newcommand{\si}{\sigma_i}
\newcommand{\la}{\lambda}
\newcommand{\eps}{\epsilon}
\newcommand{\bx}{{\bf x}}
\newcommand{\by}{{\bf y}}
\newcommand{\bJ}{{\bf J}}
\newcommand{\ts}{{\tau^{*}}}

\begin{document}

\title{Glauber Dynamics on the Cycle is Monotone}
\author{\c Serban Nacu \\ University of California, Berkeley \\ serban@stat.berkeley.edu}
\maketitle

\begin{abstract}

We study heat-bath Glauber dynamics for the ferromagnetic Ising model on a finite cycle (a graph where every vertex has degree two). We prove that the relaxation time $\tau_2$ is an increasing function of any of the couplings $J_{xy}$. We also prove some further inequalities, and obtain exact asymptotics for $\tau_2$ at low temperatures.

\end{abstract}

\section{Introduction}

Consider the ferromagnetic Ising model on a finite graph $G=(V,E)$. Let $n$ denote the number of vertices. We allow variable interaction strength, so each edge $e=[x,y] \in E$ has a coupling constant $J_{xy}>0$. As usual, every vertex $x \in V$ has a spin $\sigma_x \in \{+1,-1\}$ and the probability of a configuration is given by the Gibbs equilibrium measure $\pi(\sigma)=Z^{-1}\exp(\sum_{[x,y] \in E} J_{xy} \sigma_x \sigma_y)$. We assume no external field and no boundary conditions.

We are interested in the Glauber dynamics for this model. The Glauber dynamics is a Markov chain on the space of configurations of the Ising model, whose stationary distribution is the Gibbs measure. It comes in several flavors (heat-bath, Metropolis, etc.), but they follow the same scheme. At each step of the algorithm, one vertex is picked according to some rule (e.g. uniformly random) and with a certain probability one flips its spin (from $+1$ to $-1$ and viceversa). The flip probability depends on the value of the spin and its neighbors; if it is chosen appropiately, then the chain distribution converges to the Gibbs measure.

{}For any configuration $\sigma$ we let $\sigma^x$ be the configuration obtained by flipping spin $x$, while leaving all other spins unchanged. We shall focus on the heat-bath Glauber dynamics, in which a vertex $x$ is picked uniformly at random, then flipped with probability $\pi(\sigma^x)/(\pi(\sigma)+\pi(\sigma^x))$. Formally, this is the discrete-time Markov chain whose transition matrix $A$ has size $2^n \times 2^n$ and is defined by:

\begin{eqnarray}
A(\sigma, \sigma^x) & = & (1/n) \, \pi(\sigma^x)/(\pi(\sigma)+\pi(\sigma^x)) \nonumber \\
A(\sigma, \sigma) & = & 1 - \sum_{x \in V} A(\sigma, \sigma^x) \nonumber \\
A(\sigma, \tau) & = & 0 \quad for \, all \, other \, \tau
\end{eqnarray}

\noindent One fundamental quantity is the {\it relaxation time} $\tau_2$, defined by

\begin{equation}
\tau_2 = (n(1 - \mu_2))^{-1}
\end{equation}

\noindent where $\mu_2$ is the second largest eigenvalue of $A$. We scale by a factor of $n$ so that the $\tau_2$ defined above is the same as for a continuous-time chain where each vertex is updated at rate one. Essentially, it says we do $n$ moves of the discrete chain per unit of time.

In this paper we study the case when the graph $G$ consists of only one cycle (every vertex has degree two). Our main result is

\begin{Theorem}\label{main}
{}For the ferromagnetic Ising model on the cycle, the relaxation time $\tau_2$ is an increasing function of any of the couplings $J_{xy}$.
\end{Theorem}

This result was motivated by trying to understand how the dynamics depends on the parameters of the model. The Glauber dynamics has been used in computer simulations of the Ising model, and in describing the non-equilibrium behavior of the model. There has been a fair amount of research on its relaxation time and mixing time; see \cite{Ma} for a comprehensive survey, and \cite{KMP} for an analysis of the Glauber dynamics on trees. Most existing results give asymptotic estimates for the relaxation time when the graph size is large; in particular, there are fairly sharp estimates for trees and for the lattice ${\bf Z}^d$. However, much less is known about how the relaxation time depends on the parameters of the model, such as the temperature or the coupling constants. Such results would allow comparing mixing times for different graphs at different temperatures.

Computer simulations and heuristics suggest that, at least in the ferromagnetic case, such dependence should be monotone. Indeed, consider the dynamics started at an initial state where all spins are +1. To achieve mixing, lots of spins would have to flip to $-1$. But in the ferromagnetic model, neighboring spins tend to have the same value. The higher the interaction strength, the harder it will be for spins to flip. This led to the following

\begin{Conjecture}
(Yuval Peres, personal communication.) For the ferromagnetic Ising model on any graph, the relaxation time $\tau_2$ is an increasing function of any of the couplings $J_{xy}$.
\end{Conjecture}

What this says is that for two sets of couplings $\bJ, \bJ'$ with $J_{xy} \le J_{xy}'$, we must have $\tau_2(\bJ) \le \tau_2(\bJ')$. At first glance it may seem that one way to prove this is to find a {\it coupling} (not to be confused with the couplings $J_{xy}$) between the two Glauber dynamics, in which the dynamics for $\bJ$ is dominated by the one for $\bJ'$. This would imply the desired inequality for relaxation times.

However, such a coupling does not exist. To see why, choose some starting state (for example, all spins are +1) and let $P_n, P_n'$ be the measures obtained after $n$ steps of the dynamics. Let $\pi, \pi'$ be the corresponding Gibbs measures. We would like to define both dynamics on the same probability space so that stochastic dominance $P_n \prec P_n'$ holds. But then it must also hold in the limit, so we would need to have $\pi \prec \pi'$, which cannot hold, because both measures have mean zero. \\

In this paper, we prove the conjecture in the special case when $G$ is a cycle. The proof proceeds as follows. We show that on the cycle, the transition matrix of the Glauber dynamics preserves linear functions of the spins, and furthermore that there is a linear eigenfunction corresponding to the second largest eigenvalue. Hence by looking at the action of the dynamics on the space of linear functions of spins, the original $2^n \times 2^n$ transition matrix can be reduced to a smaller, $n\times n$ matrix, whose entries are easily computable, and mostly zeroes.

We analyze the smaller matrix to obtain a variational description of the relaxation time, as well as an inequality involving the components of its lead eigenfunction. These two combine to prove the desired monotonicity. In addition, we obtain some further inequalities involving the relaxation time, as well as low-temperature asymptotics.

\section{Two Useful Lemmas}

The results in this section are valid for the ferromagnetic Ising model on any graph. They will be needed later in the proof of our main theorem.

The transition matrix $A$ of the heat-bath Glauber dynamics was defined above. It acts on the space of functions of $\sigma$ in the natural way. It has $2^n$ eigenvalues $1=\mu_1>\mu_2 \ge \ldots \ge \mu_{2^n}$ with corresponding eigenfunctions $f_1 \equiv 1, f_2, \ldots, f_{2^n}$. Because $A$ is $\pi$-reversible we have $\pi (f_i f_j) = 0$ for $i \ne j$ (if an eigenvalue is multiple, we can choose its eigenfunctions to be $\pi$-orthogonal), and in particular $\pi f_i = 0$ for $i \ge 2$ (we use the standard notation $\pi f \equiv E_\pi f$).

An important class of functions are the {\it increasing functions}, defined in the usual way: $f$ is increasing if $f(\sigma) \le f(\tau)$ whenever $\sigma \le \tau$ (that is $\sigma_i \le \tau_i \forall i$). $f$ is strictly increasing if strict inequality holds. We have

\begin{Lemma}\label{secondeigenincr}
(Y. Peres) The second eigenvalue $\mu_2$ has an increasing eigenfunction.
\end{Lemma}

\proof The method is similar to the proof of Perron-Frobenius. We need a standard result about the Glauber dynamics: if $f$ is increasing then $Af$ is increasing (for a proof see, for example, Theorem 2.14 in~\cite{Li}). 

Let $f$ be any increasing function with $\pi f = 0$; it can be written as a linear combination of eigenvectors $f=\sum_{i=2}^{2^n} q_i f_i$. Let $\mu = \max(|\mu_2|, |\mu_{2^n}|)$. Iterating, we have 

\begin{equation}
A^m f / \mu^m = \sum_{i=2}^{2^n} q_i (\mu_i / \mu)^m f_i
\end{equation}

so if we let
\begin{eqnarray*}
g & = & \sum_{i=2}^{2^n} q_i f_i 1(\mu_i = \mu) \\
h & = & \sum_{i=2}^{2^n} q_i f_i 1(\mu_i = -\mu)
\end{eqnarray*}

\noindent we obtain convergence $A^{2m} f / \mu^{2m} \rightarrow g+h,  A^{2m+1} f / \mu^{2m+1} \rightarrow g-h$, where $g$ is some eigenvector corresponding to $\mu_2$ (a multiple of $f_2$ if $\mu_2$ is a simple eigenvalue). Since $A^nf$ is increasing for all $n$, the limits must also be increasing, so their sum must be increasing, so $g$ must be increasing. \\
It remains to show that $g$ is nonzero. If it were then we would have both $h$ and $-h$ increasing, so $h$ would also have to be zero. To avoid this, it suffices to choose $f$ not orthogonal to $f_2$ and $f_{2^n}$. For example,  $f = f_2 + f_{2^n} + K\sum \sigma_i$ has $\pi f=0$, and for large $K>0$, $\pi (f f_2) \ne 0, \pi (f f_{2^n}) \ne 0$, and $f$ is increasing. Hence $q_2 \ne 0, q_{2^n} \ne 0$ so at least one of $g$ and $h$ is nonzero.

\begin{Lemma}\label{uniqueincreigen}
(Y. Peres) If $A$ has a strictly increasing eigenfunction $f$, then $f$ corresponds to $\mu_2$.
\end{Lemma}

\proof Let $g$ be any other increasing eigenfunction. We will show there are only two possibilities: $g$ is constant (so its eigenvalue is 1), or $g$ corresponds to the same eigenvalue as $f$. Otherwise, we would have $\pi (fg)=0, \pi f=0, \pi g=0$. Since $f$ is strictly increasing, so is $f-\alpha g$ for some $\alpha >0$ small enough. Since the Ising model has the FKG property (see~\cite{FKG}), $\pi ((f - \alpha g)g) \ge \pi (f - \alpha g) \pi g$ so $\pi (fg) - \alpha \pi(g^2) \ge (\pi f)(\pi g) - \alpha (\pi g)^2$ so
$- \alpha \pi g^2 \ge 0$, contradiction. \\
Lemma~\ref{secondeigenincr} guarantees there is an increasing $g$ corresponding to $\mu_2$; hence $f$ must also correspond to $\mu_2$. \\

\section{Proof of the Main Result}

{}From now on we assume our graph $G$ is a cycle, with vertices $V=\{1, 2, \ldots, n\}$ and edges $[1, 2], [2, 3], \ldots, $ $[n, 1]$. All notation will be mod $n$, so $n+1\equiv 1$. The coupling constants are $J_{i,i+1}\equiv J_i \ge 0$. As usual, every vertex $i \in V$ has a spin $\sigma_i \in \{+1,-1\}$ and the probability of a configuration is given by $\pi(\sigma)=Z^{-1}\exp(\sum_{i=1}^n J_i \sigma_i \sigma_{i+1})$. We assume no external field and no boundary conditions.

{}For technical reasons we assume $J_i > 0$; the case when some of the $J$'s are zero can be easily obtained by limiting arguments. In particular, most results in this paper hold in the case when the graph $G$ is a straight line ($J_{n,1} = 0$).

We have already introduced increasing functions. Another important class of functions are {\it linear functions}, that is functions of the form $f(\sigma) = \sum_{i=1}^n q_i \sigma_i$ for some real constants $q_i$. On the cycle we have

\begin{Lemma}
If $f$ is linear then $Af$ is linear.
\end{Lemma}

\proof We will need the standard hyperbolic functions $\sinh(x) = 1/2(\exp(x)-\exp(-x)), \cosh(x) = 1/2(\exp(x)+\exp(-x)), \tanh(x) = \sinh(x) / \cosh(x)$. Since $A$ is a linear operator, it is enough to consider $f(\sigma) = \sigma_i$. Then $f(\s^j) = \si$ for $j \ne i$, and $f(\s^i) = -\si$, so

\begin{eqnarray*}
(Af)(\s) & = & A(\s,\s) f(\s) + \sum_j A(\s,\sigma^j) f(\sigma^j) \\
           & = & (1 - \sum_j A(\s,\sigma^j)) \si + \sum_j A(\s,\sigma^j) \si - 2A(\s, \si) \si \\
           & = & \si - 2 A(\s, \si) \si \\
           & = & (1-1/n) \si + (1/n) \si - (1/n) \frac{2\pi(\sigma^i)}{\pi(\sigma)+\pi(\sigma^i)} \si \\
           & = & (1-1/n) \si + (1/n) \si \frac{\pi(\sigma^i)-\pi(\sigma)}{\pi(\sigma^i)+\pi(\sigma)} \\
           & = & (1-1/n) \si + (1/n) \si \tanh(J_{i-1} \si \sigma_{i-1} + J_i \si \sigma_{i+1})
\end{eqnarray*}

Since $\tanh$ is an odd function, the second term does not depend on $\si$, so

\begin{equation}
(Af)(\s) = (1-1/n) \si + (1/n) \tanh(J_{i-1} \sigma_{i-1} + J_{i+1} \sigma_{i+1})
\end{equation}

and the $\tanh$ is easily expressed as a linear function 

\begin{equation}
\tanh(J_{i-1} \sigma_{i-1} + J_{i+1} \sigma_{i+1}) = a\sigma_{i-1} + b\sigma_{i+1}
\end{equation}

since $\sigma_{i-1}$ and $\sigma_{i+1}$ only take the values $+1$ and $-1$, and both sides of the equality are odd functions.
So we only need to check that $a+b=\tanh(J_{i-1}+J_i), a-b=\tanh(J_{i-1}-J_i)$.
After a bit of algebra, this yields 

\begin{eqnarray}
a & = & \sinh(2J_{i-1}) / (\cosh(2J_{i-1}) + \cosh(2J_i)) \nonumber \\
b & = & \sinh(2J_i) / (\cosh(2J_{i-1}) + \cosh(2J_i)) \nonumber
\end{eqnarray}

\noindent and the lemma is proven. \\

Hence we can look at the restriction of A to the space of linear functions. With respect to the basis $\{\sigma_1, \sigma_2, \ldots, \sigma_n\}$, this is described by an $n \times n$ matrix $L$ whose entries are all zero, except on the diagonal, and immediately above and below it. Let $s_i = \sinh(2J_i)$ and $c_i = \cosh(2J_i)$. Then we have 

\begin{eqnarray*}
L(i,i) &=& 1-1/n \\ L(i,i-1) &=& (1/n)s_{i-1} / (c_{i-1} + c_i) \\ L(i,i+1) &=& (1/n)s_i / (c_{i-1} + c_i)
\end{eqnarray*}

\noindent Note that an eigenvalue of $L$ is also an eigenvalue of $A$, so all must be real. Since all entries on the diagonal of $L$ are equal, we can write $L=(1-1/n)I + (1/n)M$, where the diagonal of $M$ is all zeroes, and the only non-zero entries are

\begin{equation}
M(i,i-1) = s_{i-1} / (c_{i-1} + c_i), \quad M(i,i+1) = s_i / (c_{i-1} + c_i)
\end{equation}

\noindent Let the eigenvalues of $M$ be $\la_1 \ge \ldots \ge \la_n$. Then the eigenvalues of $L$ are $1-(1-\la_i)/n$ for $1 \le i \le n$, and $M$ and $L$ have the same eigenvectors. It is easy to see $M$ is irreducible. Hence by Perron-Frobenius, $\la_1 > 0$, $\la_1$ is a simple eigenvalue, and there exists a {\it dominant eigenvector} $\bx$ (that is, corresponding to $\la_1$) which has non-negative entries. In fact, the entries must be strictly positive, since we have $M\bx = \la_1 \bx$, so

\begin{equation}
M(i,i-1)x_{i-1}+M(i,i+1)x_{i+1} = \la_1 x_i
\end{equation}

\noindent so if $x_i=0$ for some $i$, then also $x_{i+1}=0$, so $\bx = 0$, contradiction. But then $f(\s) = \sum x_i \si$ is a strictly increasing function of $\sigma$, and also an eigenfunction of $A$. Lemma~\ref{uniqueincreigen} then implies that $f$ corresponds to the second eigenvalue. We have obtained

\begin{Corollary}
$\mu_2 = 1 - (1 - \la_1) / n$, and $\mu_2$ has a linear eigenfunction.
\end{Corollary}

Hence the relaxation time of the Glauber dynamics can be computed as $\tau_2 = 1 / (n(1 - \mu_2)) = 1 / (1 - \la_1)$, so from now on we can focus on the small matrix $M$. We will need an inequality which roughly says that its lead eigenvectors are not too far from a constant vector.

\begin{Lemma}\label{mainineq}
If $\bx$ is a non-negative dominant eigenvector of $M$, then
\begin{equation}\label{eigenvectorineq} 
(c_i - 1) / s_i \le x_i / x_{i+1} \le (c_i + 1) / s_i
\end{equation}
\end{Lemma}

\proof The method could be called ``cyclic induction.'' Note that $c_i^2 - s_i^2 = 1$, so $(c_i - 1) / s_i = s_i / (c_i + 1) < 1$. We have $M\bx = \la_1 \bx$, so

\begin{equation}\label{eigenvectoreq}
\la_1 (c_{i-1} + c_i) x_i = s_{i-1} x_{i-1} + s_i x_{i+1} \quad \forall i
\end{equation}

Assume the left-hand inequality of~(\ref{eigenvectorineq}) holds for $i-1$ so $(c_{i-1} - 1) / s_{i-1} \le x_{i-1} / x_i$. Combining this with (\ref{eigenvectoreq}), we get (since $\la_1 < 1$)

\begin{displaymath}
(c_{i-1} + c_i) x_i \ge \la_1 (c_{i-1} + c_i) x_i \ge (c_{i-1} - 1) x_i + s_i x_{i+1}
\end{displaymath}

\noindent so $(c_i + 1) x_i \ge s_i x_{i+1}$ so $x_i / x_{i+1} \ge s_i / (c_i + 1) = (c_i - 1) / s_i$. Hence the inequality holds for $i$ as well; hence it holds either for all $i$ or for none. But in the latter case we would get $x_i / x_{i+1} < (c_i - 1) / s_i < 1$ for all $i$, hence $x_1 < x_2 < \ldots < x_n < x_1$, which is impossible.

The right-hand inequality can be proven in a similar way: if the inequality $(c_i + 1) / s_i \le x_i / x_{i+1}$ holds for $i-1$, it also holds for $i$; it cannot hold for all $i$, so it cannot hold for any $i$. Alternatively, the right-hand side can be deduced from the left-hand side by symmetry. \\

An immediate consequence of this is
\begin{Proposition}\label{mainineq2}
If $\bx$ is a non-negative dominant eigenvector of $M$, then $x_i/x_{i+1} + x_{i+1}/x_i \le 2c_i/s_i$.
\end{Proposition}

\proof The function $f(t) = t + 1/t$ is decreasing for $t<1$ and increasing for $t>1$, so it follows from Lemma~\ref{mainineq} that $f(x_i/x_{i+1}) \le f((c_i+1)/s_i)$ if $x_i \ge x_{i+1}$, and
$f(x_i/x_{i+1}) \le f((c_i-1)/s_i)$ if $x_i \le x_{i+1}$. Both upper bounds are equal to $2c_i/s_i$. \\

To apply the last result we need an expression for the largest eigenvalue in terms of its eigenvector. Note that the product of the elements under the diagonal of $M$ is equal to the product of the ones above the diagonal, so $M$ is reversible in the following sense. Define $q(i) = c_i + c_{i+1}$. Then $q(i)M(i,j)=q(j)M(j,i) \forall i,j$. So if $Q$ is the diagonal matrix with $Q(i,i)=q(i)$ then $Q^{1/2}MQ^{-1/2}$ is symmetric, and has the same eigenvalues as $M$. Then by the usual extremal characterization,

\begin{displaymath}
\la_1 = \sup \frac{(Q^{1/2}MQ^{-1/2}v,v)}{(v,v)}
\end{displaymath}

\noindent where the $\sup$ can be taken over all nonzero vectors $v$, and is attained at dominant eigenvectors of $Q^{1/2}MQ^{-1/2}$. The numerator expands as 

\begin{eqnarray*}
(Q^{1/2}MQ^{-1/2}v,v) & = & \sum 2v_i v_{i+1} \sqrt{M(i,i+1)M(i+1,i)} \\
                              & = & \sum 2v_i v_{i+1} s_i / \sqrt{(c_{i-1}+c_i)(c_i + c_{i+1})}
\end{eqnarray*}

\noindent so if we finally substitute $v_i = (c_i + c_{i+1})^{1/2} y_i$ we obtain

\begin{Lemma}\label{extrchar}
$\la_1 = \sup_{\bf y} \, (\sum 2y_i y_{i+1} s_i) / (\sum (y_i^2 + y_{i+1}^2) c_i)$ and the $\sup$ is attained at dominant eigenvectors of $M$.
\end{Lemma}

\noindent Now let $h(\by, \bJ) = (\sum 2y_i y_{i+1} s_i) / (\sum (y_i^2 + y_{i+1}^2) c_i)$ and let $\bx \equiv \bx(\bJ)$ be a non-negative dominant eigenvector of $M$, chosen so $\bx(\cdot)$ is a smooth function in a neighborhood of $\bJ$. This can be done because $\la_1$ is a simple eigenvalue; see, for example, \cite{Ka}, Chapter 2.1.
We have $\la_1(\bJ) = h(\bx(\bJ), \bJ)$. Because for fixed $\bJ$, $\sup h(\by, \bJ)$ is attained at $\by = \bx(\bJ)$, we have

\begin{equation}
\partial h / \partial y_k (\bx(\bJ), \bJ) = 0 \quad \forall k
\end{equation}

so by the chain rule
\begin{equation}
\partial \la_1 / \partial J_i = \partial h / \partial J_i + \sum_k (\partial h / \partial y_k) (\partial x_k / \partial J_i) = \partial h / \partial J_i 
\end{equation}

Hence

\begin{eqnarray}
\partial \la_1 / \partial J_i & = & \partial h / \partial J_i (\bx, \bJ) \nonumber \\
                                     & =  &  2(2x_i x_{i+1} c_i - (x_i^2 + x_{i+1}^2) s_i h(\bx, \bJ) ) / \sum (x_i^2 + x_{i+1}^2) c_i \label{derivative}
\end{eqnarray}

{}From Proposition~\ref{mainineq2}, $(x_i^2 + x_{i+1}^2) / (2x_i x_{i+1}) \le c_i / s_i$ and since $h(\bx, \bJ) = \la_1 < 1$ it follows that $\partial \la_1 / \partial J_i > 0$. This completes the proof of Theorem~\ref{main}.

\section{Further Inequalities}

We can get more precise bounds on the derivatives of $\la_1$. Indeed, equation (\ref{derivative}) and Proposition~\ref{mainineq2} yield easily
\begin{equation}
\partial \la_1 / \partial J_i \ge 2(1 - \la_1) (2x_i x_{i+1} c_i) / \sum (x_i^2 + x_{i+1}^2) c_i
\end{equation}

hence 
\begin{eqnarray}
\sum (s_i / c_i) \partial \la_1 / \partial J_i  & \ge  & 2(1 - \la_1) \sum (2x_i x_{i+1} s_i) / \sum (x_i^2 + x_{i+1}^2) c_i \nonumber \\
& = & 2(1- \la_1) \la_1
\end{eqnarray}

So we can obtain an interesting inequality by differentiating along a particular direction. To do that, consider functions $J_i(t)$ with $J_i(0) = J_i$ and $J_i'(t) = \tanh(2J_i(t))$. This differential equation is easy to solve and yields 

\begin{equation}
J_i(t)  = (1/2) \sinh^{-1} (\exp(2t) \sinh(2J_i))
\end{equation}

\noindent (here $\sinh^{-1}$ denotes the inverse hyperbolic function). Then the function $\la(t) = \la_1(J_1(t), \ldots, J_n(t))$ satisfies $\la' \ge 2\la(1 - \la)$. This becomes even simpler if we substitute 

\begin{equation}
\ts = \tau_2 - 1 = \la / (1 - \la)
\end{equation}

We obtain $\ts ' \ge 2\ts$ so $(\log \ts)' \ge 2$ and $\ts (t) \ge \exp(2t)$. Hence, with a slight abuse of notation (by $f(\bJ)$ we mean that $f$ is applied to all entries of $\bJ$), we have

\begin{Proposition}
$ \ts((1/2)\sinh^{-1} (\exp(2t) \sinh(2\bJ))) \ge \exp(2t) \ts(\bJ) \forall \bJ, t \ge 0$.
\end{Proposition}

{}For large $x$, $\sinh^{-1}(x)$ is close to $\log(x)$, so we expect $J_i(t)$ to be approximately $J_i + t$. In fact it is easy to prove that $\exp(2t) \sinh(2u) \le \sinh(2t+2u)$, so $J_i(t) \le J_i + t$,
so from Theorem~\ref{main} we get

\begin{Theorem}
$\ts(\bJ + t) \ge \exp(2t) \ts(\bJ)$ for all $t \ge 0$.
\end{Theorem}

As above, $\bJ + t$ means $t$ is added to all entries of $\bJ$. Numerical computations confirm that the factor 2 in the exponent is optimal. It would be nice to have a similar inequality where the elements of $\bJ$ are {\it multiplied} by $t$, so that we can estimate the effect of a change of temperature on the relaxation time. Note, however, that $t\bJ > t + \bJ$ for large $t$, so the inequality we have proven is in some sense stronger.

\section{Low Temperature Asymptotics}

Introduce the inverse temperature $\beta$; the coupling constants are now $\beta J_1, \ldots, \beta J_n$. We are interested in $\tau_2(\beta \bJ)$ for large $\beta$ and fixed $\bJ$. This turns out to be asymptotically exponential in $\beta$. More precisely, let $J_{max} = \max J_i$ and $J_{min} = \min J_i$. We have

\begin{Proposition}
The following holds as $\beta \rightarrow \infty$:
$$\lim \tau_2(\beta \bJ) / \exp(C\beta) = D$$
\noindent where $C = 2(J_{max}+J_{min})$ and $D=(1/2) \# \{i:J_i=J_{max}\} / \# \{i:J_i = J_{min} \} $.
\end{Proposition}

\proof As before, let $\bx$ be a dominant eigenvector. Lemma~\ref{mainineq} guarantees that the ratios $x_i / x_{i+1}$ are close to 1 for large $\beta$. That is, given $\eps > 0$, $1 - \eps \le x_i / x_{i+1} \le 1 + \eps$ for $\beta$ large enough. We can set $x_1 = 1$ so we get $(1 - \eps)^n \le x_i \le (1 + \eps)^n \, \forall i$.

If we apply Lemma~\ref{extrchar} for $\by \equiv 1$ we obtain $\la_1 \ge \sum s_i / \sum c_i$ so
\begin{equation}
\tau_2 = 1 / (1 - \la_1) \ge \sum c_i / \sum (c_i - s_i)
\end{equation}

If we apply Lemma~\ref{extrchar} for $\by = \bx$ and substitute $\tau_2 = 1 / (1 - \la_1)$, we obtain
\begin{eqnarray}
\tau_2 &=& \sum c_i(x_i^2 + x_{i+1}^2) / [ \sum c_i(x_i - x_{i+1})^2 + \sum 2x_i x_{i+1} (c_i - s_i) ] \nonumber \\
          & \le & ((1+\eps)/(1-\eps))^{2n} \sum c_i / \sum (c_i - s_i)
\end{eqnarray}

\noindent The ratio $\sum c_i / \sum (c_i - s_i)$ is easily shown to be asymptotically $D \exp(C \beta)$. Since $\eps$ is arbitrary, the proof is complete. \\

\noindent {\bf Acknowledgments.} Many thanks to Yuval Peres, my thesis advisor. His guidance was instrumental in obtaining these results. I also thank him for permission to include Lemmas~\ref{secondeigenincr} and~\ref{uniqueincreigen}. Thanks also to Noam Berger, 
Elchanan Mossel and Dror Weitz for their comments. \\

\end{document}